\documentclass{amsart}

\usepackage{amssymb,amsfonts,latexsym}

\numberwithin{equation}{section}

\newtheorem{theorem}{Theorem}[section]

\newtheorem{lemma}[theorem]{Lemma}

\newtheorem{proposition}[theorem]{Proposition}

\newtheorem{corollary}[theorem]{Corollary}

\theoremstyle{definition}

\newtheorem{remark}[theorem]{Remark}

\newtheorem *{Theorem A}{Theorem A}
\newtheorem *{Corollary B}{Corollary B}

\newcommand{\kb}{{K_{\beta}}}

\newcommand{\ach}{\check{A}}

\newcommand{\la}{\langle\,}

\newcommand{\ra}{\,\rangle}

\newcommand{\ben}{\begin{enumerate}}

\newcommand{\een}{\end{enumerate}}

\newcommand{\im}{{\rm im}}

\newcommand{\res}{{\rm res}}

\newcommand{\infl}{{\rm inf}}

\newcommand{\C}{{\mathbb C}}

\begin{document}

\title
[on groups of central type]
{On groups of central type, non-degenerate and bijective cohomology classes}

\author{Nir Ben David}

\address{Department of Mathematics, Technion-Israel Institute of Technology, Haifa 32000, Israel}
\email{benda@tx.technion.ac.il}

\author{Yuval Ginosar}

\address{Department of Mathematics, University of Haifa, Haifa 31905, Israel}
\email{ginosar@math.haifa.ac.il}

\date{\today}


\begin{abstract}
A finite group $G$ is of central type (in the non-classical sense)
if it admits a non-degenerate cohomology class $[c]\in H^2(G,\C^*)$ ($G$ acts trivially on $\C^*$).
Groups of central type play a fundamental role in the classification of
semisimple triangular complex Hopf algebras and can be determined by their
repre\-sen\-ta\-tion-\-theo\-retical properties.

Suppose that a finite group $Q$ acts on an abelian group $A$ so that there exists a
bijective 1-cocycle $\pi\in Z^1(Q,\ach)$, where $\ach=\rm{Hom}(A,\C^*)$ is
endowed with the diagonal $Q$-action.
Under this assumption,
Etingof and Gelaki gave an explicit formula for a non-degenerate 2-cocycle in $Z^2(G,\C^*)$,
where $G:=A\rtimes Q$.
Hence, the semidirect product $G$ is of central type.

In this paper we present a more general correspondence
between bijective and non-degenerate cohomology classes. In particular, given a
bijective class $[\pi]\in H^1(Q,\ach)$ as above, we construct non-degenerate classes
$[c_{\pi}]\in H^2(G,\C^*)$ for certain extensions
$1\to A\to G\to Q\to 1$ which are not necessarily split.
We thus strictly extend the above family of central type groups.
\end{abstract}

\maketitle

\section{Introduction}

It is well known that the
dimension of an irreducible complex representation of a finite group $\Gamma$ cannot
exceed $\sqrt{[\Gamma:Z(\Gamma)]}$, where $Z(\Gamma)$ is the
center of $\Gamma$ (see e.g. [I, Corollary 2.30]).
Classically, $\Gamma$ is termed ``of central
type" if it admits an irreducible representation $V$ of dimension $n:=\sqrt{[\Gamma:Z(\Gamma)]}$.
However, here (following [AHN]) we call the quotient group $G:=\Gamma/Z(\Gamma)$
a group of {\em central type}.

Let $V$ be an $n$-dimensional irreducible representation of $\Gamma$.
Since $Z(\Gamma)$ is mapped to the scalar matrices, the representation $V$ determines
a projective representation $\bar{V}$ of $G$.
Let $[c']\in H^2(G,Z(\Gamma))$ be the cohomology class corresponding to the
extension $1\to Z(\Gamma)\to \Gamma\to G\to 1$, and let
$[c]\in H^2(G,\C^*)$  be the image of $[c']$ under the induced
map $H^2(G,Z(\Gamma))\to H^2(G,\C^*)$ (where $G$ acts trivially on $\C^*$).
Then $[c]$ is associated to the projective representation $\bar{V}$
and the twisted group algebra $\C^{c}[G]$ is isomorphic to $M_n(\C)$.
It turns out [DMJ, Corollary 3] that the action of $G$ on $\C^{c}[G]$ by conjugation is isomorphic
to the regular representation of $G$. This fact will be useful in the sequel.
A 2-cocycle satisfying the above properties is called {\it non-degenerate}.

A remarkable result of Howlett and Isaacs asserts (assuming the
classification of finite simple groups) that a group of central
type $G$ is solvable [HI, Theorem 7.3].

Groups of central type have attracted considerable attention due to their
represen\-tation-\-theoretical properties and their fundamental role in
the theory of fusion symmetric categories via the theory of
semisimple triangular Hopf algebras. More specifically, it was
proved that the isomorphism classes of semisimple
triangular Hopf algebras of dimension $m$ over $\C$ are in
bijection with the isomorphism classes of quadruples
$(\Gamma,H,c,u)$ where $\Gamma$ is a group of order $m$, $H$ is a central type
subgroup of $\Gamma$, $c$ is a non-degenerate $2$-cocycle on $H$ and $u$
in an involution in $\Gamma$ [EG1,3].
The above bijection is given by twists for finite groups in the sense of Drinfeld [D].
Such a twist is supported on a subgroup of central type,
and any finite dimensional semisimple triangular Hopf algebra over $\C$
is obtained from a group of central type equipped with a twist [EG3].

This motivates our interest in the explicit construction of central type
groups. An important family of such groups was described in [EG2,3]: given a finite group $Q$
acting on an abelian group $A$, $\ach=\rm{Hom}(A,\C^*)$ is a $Q$-module
under the diagonal action (see (\ref{pairact})).
Let $\pi:Q\to \ach$ be a bijective 1-cocycle, so that $\pi(g_1g_2)=\pi(g_1)g_1(\pi(g_2))$
for every $g_1,g_2\in Q$, and $|Q|=|A|$. Etingof and Gelaki then construct a non-degenerate
2-cocycle $c$ on the semidirect product $G:=A\rtimes Q$, showing that $G$ is a group of central type.
This non-degenerate 2-cocycle has an extra property: its restriction to $A$ is trivial.

Our main goal in this paper is to generalize the construction in [EG3]
to extensions $1\to A\to G\to Q\to 1$ which are not necessarily split.
We also place the correspondence between bijective 1-cocycles $\pi: Q\to\ach$ and
non-degenerate 2-cocycles on $G$ in a more general context, as follows.

Let
\begin{equation}\label{exten}
[\beta]:1\to A\to G\to Q\to 1, \ \ [\beta]\in H^2(Q,A)
\end{equation}
be an extension of a finite group $Q$ by a finite abelian group $A$. Let
$$
\res^G_A:H^2(G,\C^*)\to H^2(A,\C^*),\ \ \infl^Q_G:H^2(Q,\C^*)\to H^2(G,\C^*)
$$
be the restriction
and inflation maps respectively, where $G$ acts trivially on $\C^*.$

As noted in the beginning of Section \ref{2main}, if a 1-cocycle is bijective, so are all
elements in its cohomology class; similarly for a non-degenerate 2-cocycle. Moreover,
if $|A|=|Q|$ and $c:G\times G \to \C^*$ is a non-degenerate 2-cocycle, then so is
$c\cdot f$ for every 2-cocycle $f$ inflated from $Q$. From now on we shall therefore refer
also to cohomology classes (rather than cocycles) which are bijective or non-degenerate
modulo the image of $\inf^Q_G$.

Our main result is:
\begin{Theorem A}\label{main2}
Let (\ref{exten}) be an extension of finite groups such that $|A|=|Q|$.
Then there is a 1-1 correspondence between bijective classes $[\pi]\in H^1(Q,\ach)$
such that $[\beta] \cup [\pi]=0\in H^3(Q,\C^*)$ and
non-degenerate classes in $\ker(\res^G_A)\mod[\im(\inf^Q_G)]$.
\end{Theorem A}

When $G=A\rtimes Q$, then clearly any $[\pi]\in H^1(Q,\ach)$ satisfies $[\beta] \cup [\pi]=0$,
and the correspondence in Theorem A yields the cohomology class of the non-degenerate
2-cocycle that was constructed in [EG3].

Theorem A shows that if $A$ is a normal abelian subgroup of $G$ such that
$|A|=|G/A|$, then any non-degenerate class in $\ker(\res^G_A)$ is obtained,
up to inflation, from a certain bijective class.

Further, fix any bijective class $[\pi]\in H^1(Q,\ach)$.
The correspondence in Theorem A gives rise to non-degenerate cohomology classes
in $H^2(G,\C^*)$ for all extensions (\ref{exten}) such that $[\beta] \cup [\pi]=0$.
We conclude:

\begin{Corollary B}\label{cor}
Let $A$ be a finite abelian group, $Q$ a finite group acting on
$A$ and $[\pi]\in H^1(Q,\ach)$ a bijective class.
Then for every $[\beta]\in H^2(Q,A)$ such that $[\beta] \cup [\pi]=0$,
the group $G$ determined by the extension $[\beta]: 1\to A\to G\to Q\to 1$
is of central type.
\end{Corollary B}

Corollary B gives a way to construct groups of central type from bijective cohomology classes.
In practice, it yields groups of central type that cannot be constructed as a semidirect product in any way.
In Section \ref{Example} we give an example of a non-split extension
$[\beta]: 1\to A\to G \to Q\to 1$ and a bijective class $[\pi] \in H^1(Q,\check{A})$ such that $[\beta] \cup [\pi]=0$.
The group $G$ (of order $64$) is hence of central type, but does not contain any abelian
normal subgroup $N$ of order $|G/N|$ such that $G:=N\rtimes (G/N)$.

The 1-1 correspondence in Theorem A is obtained from an isomorphism of the quotient $\ker(\res^G_A)/\im(\inf^Q_G)$ and
$$\kb:=\{[\pi]\in H^1(Q,\check{A})\ \ |\ \ [\beta] \cup [\pi]=0\in H^3(Q,\C^*)\} .$$
There is a standard identification of both $\ker(\res^G_A)/\im(\inf^Q_G)$ and $\kb$
with a limit term of the Lyndon-Hochschild-Serre (LHS) spectral sequence (see \S \ref{recall} hereafter).
We shall however give a self contained description of mutually inverse isomorphisms
\begin{equation} \label{pic}
\Pi:\ker(\res^G_A)/\im(\infl^Q_G)\to \kb,\ \ C:\kb \to \ker(\res^G_A)/\im(\infl^Q_G)
\end{equation}
in terms of representatives and cocycles. This description will be convenient for
establishing the correspondence in Theorem A.

The paper is organized as follows:

After some cohomological background in Section \ref{pre},
we describe the mutually inverse isomorphisms $\Pi$ and $C$ in Section \ref{1main}.
These isomorphisms induce the 1-1 correspondence in Theorem A as shown in Section \ref{2main}.
The isomorphism $\Pi$, essentially that of \cite[Section 1.7]{K},
is interpreted here via the action of $Q$ on the primitive idempotents of the group algebra
$\C[A]$.  Section \ref{Example} is the example mentioned above.
Finally, in Section \ref{twist} we apply the 1-1 correspondence between bijective and
non-degenerate cohomology classes in Theorem A to construct twists for finite groups explicitly.\\

{\bf Acknowledgments.}
Part of this work is based on the M.Sc.~thesis of the first author
under the supervision of E.~Aljadeff and S.~Gelaki.
We are grateful to them for their significant contribution.
We also thank D.~Blanc, H.W.~Henn, M.~Natapov,
\'{A}.~del R\'{\i}o and U.~Weiss for their help.\\
The first author was supported by Technion V.P.R. Fund - Dent Charitable Trust
and by The Israel Science Foundation (grant No. 70/02-1).

\section{Preliminaries}\label{pre}

\subsection{}
We first recall some elementary notation.
Let $Q$ be a finite group acting on a finite abelian group $A$. This action induces a diagonal
action of $Q$ on $\ach=\rm{Hom}(A,\C^*)$ defined by the pairing between $\ach$ and $A$:
\begin{equation}\label{pairact}
\la g(\chi),a\ra =\la \chi,g^{-1}(a)\ra \in \C^*,
\end{equation}
where $g\in Q, \chi\in \ach$ and $a\in A$ (the action on $\C^*$ is trivial).

Let (\ref{exten})
be an extension of $Q$ by $A$.
Let $\{\bar{g}\}_{g\in Q}$ be a transversal set for $Q$ in $G$.
Then any element in $G$ is uniquely expressed as $a\bar{g}$,
where $a\in A$ and $g\in Q$. Recall that the multiplication in $G$
is determined by the conjugation of elements in $A$ via the given action
$$\bar{g}a\bar{g}^{-1}=g(a), \ \ a\in A, g\in Q,$$
and by the 2-cocycle $\beta:Q\times Q\to A$ as follows
\begin{equation}\label{beta}
\bar{g}_1\bar{g}_2=\beta(g_1,g_2)\overline{g_1g_2},\ \ g_1,g_2\in Q.
\end{equation}
We may further assume that $\beta$ is normalized, i.e.
$$
\beta(g,1)=\beta(1,g)=1,\ \ \forall g\in Q,
$$
by choosing the trivial element of $G$ as a representative of the trivial element in $Q$. \\

\subsection{} \label{recall}

The second cohomology group $H^2(G,\ach)$ admits the filtration
\begin{equation}
\im(\infl^Q_G)\subseteq \ker(\res^G_A) \subseteq H^2(G,\ach),
\end{equation}
and it is well known that the quotient
$\ker(\res^G_A)/\im(\inf^Q_G)$ is isomorphic to the term
$E^{1,1}_{\infty}$ in the LHS spectral sequence for the extension
(\ref{exten}). Moreover, $E^{1,1}_{\infty}$ canonically embeds into
$E^{1,1}_2$ as the kernel of the differential $d^{1,1}_2:E^{1,1}_2
\rightarrow E^{3,0}_2$. Since $A$ is abelian and acts trivially on
$\C^*$, we have
$$E^{1,1}_2=H^1(Q,H^1(A,\C^*))=H^1(Q,\ach), \ \ E^{3,0}_2=H^3(Q,H^0(A,\C^*))=H^3(Q,\C^*).$$
Furthermore, by [HS, Theorem 4], the differential $d^{1,1}_2$ in this case amounts
to multiplication by the class of the extension up to sign. More precisely,
\begin{eqnarray}
\begin{array}{rrcl}
 d^{1,1}_2: & H^1(Q,\ach) & \to &  H^3(Q,\C^*)\\
  & [\pi] & \mapsto & -[\beta] \cup [\pi].
\end{array}
\end{eqnarray}
where $[\beta]\in H^2(Q,A)$ is determined by the extension and
$\cup$ is the usual cup product followed by the pairing $A\otimes\ach\to \C$ (see (\ref{cup})).
Both $\kb$ and $\ker(\res^G_A)/\im(\inf^Q_G)$ may therefore be identified with $E^{1,1}_{\infty}$.

A detailed study of $d^{1,1}_2$, as well as of other low degree differentials
in the general setup can be found in [H].

\section{The Isomorphisms $C$ and $\Pi$} \label{1main}

In this section we explicitly describe the mutually inverse isomorphisms
$C:\kb\to\ker(\res^G_A)/\im(\inf^Q_G)$ (\S \ref{subsecC}) and
$\Pi:\ker(\res^G_A)/\im(\inf^Q_G)\to \kb$ (\S \ref{subsecPi}).

\subsection{}\label{subsecC}
We first construct the isomorphism $C:\kb\to\ker(\res^G_A)/\im(\inf^Q_G)$.
Let $\pi\in Z^1(Q,\ach)$ be any 1-cocycle. Then it determines a 2-cochain $\varphi_{\pi}\in C^2(G,\C^*)$
as follows (compare with the construction for semidirect products [EG3, section 8]).
For $\gamma_1=a_1\bar{g}_1,\gamma_2=a_2\bar{g}_2\in G$, let
\begin{equation}
\varphi_{\pi}(\gamma_1,\gamma_2):=\la \pi(g_1),g_1(a_2)\ra ^{-1}\ \in \C^*.
\end{equation}
Equivalently, by (\ref{pairact}) and the 1-cocycle condition we obtain
\begin{equation}\label{obt}
\varphi_{\pi}(\gamma_1,\gamma_2)=\la \pi(g_1^{-1}),a_2\ra .
\end{equation}
Note that for every $\pi_1,\pi_2\in Z^1(Q,\ach)$
\begin{equation}\label{varphimulti}
\varphi_{\pi_1\pi_2}=\varphi_{\pi_1}\varphi_{\pi_2}.
\end{equation}
In general, the cochain $\varphi_{\pi}$ does not satisfy the
cocycle condition. However, when $[\pi]$ is in $\kb\subseteq H^1(Q,\ach)$, we can modify its construction to
make it a cocycle, as follows.
Recall that for any (left) $G$-module $M$, the $n$-th coboundary map
$\delta_G^n:C^n(G,M)\to C^{n+1}(G,M)$ is defined by
\begin{equation*}
\begin{split}
\delta_G^n(\varphi)(\gamma_1,&\dotsc,\gamma_{n+1}):=~
~\gamma_1\cdot\varphi(\gamma_2\dotsc,\gamma_{n+1})\\
&+~\sum\limits_{i=1}^{n}(-1)^i\varphi(\gamma_1,\dotsc,\gamma_i\gamma_{i+1},\dotsc,\gamma_{n+1})
+~(-1)^{n+1}\varphi(\gamma_1,\dotsc,\gamma_{n}),
\end{split}
\end{equation*}
for every $\varphi$ in $C^n(G,M)$.

\begin{proposition} \label{vp}
For any $\pi\in Z^1(Q,\ach)$ and any $a_1\bar{g}_1,a_2\bar{g}_2,a_3\bar{g}_3\in G$,\\
$\delta_G^2(\varphi_{\pi})(a_1\bar{g}_1,a_2\bar{g}_2,a_3\bar{g}_3)=
\la \pi(g_1),g_1(\beta(g_2,g_3))\ra ^{-1}.$
\end{proposition}
\begin{proof}
This follows directly from the definition of $\delta_G^2$.
\end{proof}
Next, let $[\pi]\in \kb$. Then $[\beta] \cup [\pi]=0\in H^3(Q,\C^*)$ -
i.e., the 3-cochain
\begin{equation} \label{cup}
(g_1,g_2,g_3)\mapsto\la \pi(g_1),g_1(\beta(g_2,g_3))\ra
\end{equation}
is a 3-coboundary on $Q$ [HS, p.118].
Hence there exists a 2-cochain $\zeta_{\pi}\in C^2(Q,\C^*)$ such that
\begin{equation}  \label{vs}
\delta_Q^2(\zeta_{\pi})(g_1,g_2,g_3)=\la \pi(g_1),g_1(\beta(g_2,g_3))\ra ,
\forall g_1,g_2,g_3\in Q.
\end{equation}
We now construct a 2-cocycle on $G$ from the cochain $\varphi_{\pi}$.
Let
\begin{equation}\label{cp}
c_{\pi}:=\varphi_{\pi}\bar{\zeta}_{\pi},
\end{equation}
where $\bar{\zeta}_{\pi}\in C^2(G,\C^*)$ is the inflation of
$\zeta_{\pi}$ to $G$. By Proposition \ref{vp} and
(\ref{vs}), we have
$\delta_G^2(\varphi_{\pi}\bar{\zeta}_{\pi})=0$. Consequently
\begin{equation}
c_{\pi}\in Z^2(G,\C^*).
\end{equation}
We need the following
\begin{proposition}\label{condpi}
Let $c_{\pi}$ be as above. Then:
\begin{enumerate}
\item The restriction of $c_{\pi}$ to $A$ is cohomologically trivial.
\item Any other choice of a 2-cochain $\zeta_{\pi}'\in C^2(Q,\C^*)$
(satisfying (\ref{vs}))
yields a 2-cocycle $c_{\pi}'\in Z^2(G,\C^*)$, with $c_{\pi}'c_{\pi}^{-1}$ inflated from $Q$.
\item Any other choice of a transversal set $\{\bar{g}'\}_{g\in Q}$ for $Q$ in $G$ yields a
2-cocycle $c_{\pi}'\in Z^2(G,\C^*)$, with $c_{\pi}'c_{\pi}^{-1}$ inflated from $Q$.
\item  If $\pi\in B^1(Q,\ach)$ then $c_{\pi}\in B^2(G,\C^*)$.
\item For every $\pi_1,\pi_2\in Z^1(Q,\ach)$, the 2-cocycle $c_{\pi_1}c_{\pi_2}c_{\pi_1\pi_2}^{-1}$
is inflated from $Q$.
\end{enumerate}
\end{proposition}
\begin{proof}
(1) By (\ref{obt}), for every $a_1,a_2\in A$,
$$
\varphi_{\pi}(a_1,a_2)=\la \pi(1),a_2\ra =1.
$$
Next, since the cochain $\bar{\zeta}_{\pi}$ is inflated from $Q$,
its restriction to $A$ is a constant 2-cocycle, and in particular a coboundary.
It follows that the restriction of $c_{\pi}=\varphi_{\pi}\bar{\zeta}_{\pi}$
to $A$ is a coboundary.\\
\noindent(2) This holds since $\zeta_{\pi}'$ differs from $\zeta_{\pi}$ by a 2-cocycle on $Q$. \\
\noindent(3) The two transversal sets differ by a 1-cochain $\lambda\in C^1(Q,A)$, i.e.
$\bar{g}=\lambda(g)\bar{g}'$ for every $g\in Q$.
Let
$$
\gamma_1=a_1\bar{g_1}=a_1\lambda{(g_1)}\bar{g_1}',
\gamma_2=a_2\bar{g_2}=a_2\lambda{(g_2)}\bar{g_2}' \in G.
$$
Then
$$
\begin{array}{l}
c'_{\pi}(\gamma_1,\gamma_2)=(\varphi'_{\pi}\bar{\zeta_{\pi}}')(\gamma_1,\gamma_2)=
\la \pi(g_1^{-1}),a_2\lambda{(g_2)}\ra \cdot\zeta'_{\pi}(g_1,g_2)=\\
=\la \pi(g_1^{-1}),a_2\ra \cdot\la \pi(g_1^{-1}),\lambda{(g_2)}\ra \cdot\zeta'_{\pi}(g_1,g_2)=\\
(\varphi_{\pi}\bar{\zeta_{\pi}})(\gamma_1,\gamma_2)\cdot\zeta_{\pi}(g_1,g_2)^{-1}\cdot
\la \pi(g_1^{-1}),\lambda{(g_2)}\ra \cdot\zeta'_{\pi}(g_1,g_2)=
c_{\pi}(\gamma_1,\gamma_2)\eta_{\pi}(\gamma_1,\gamma_2),
\end{array}
$$
where $\eta_{\pi}\in Z^2(G,\C^*)$ defined by
$$
\eta_{\pi}(\gamma_1,\gamma_2):=\zeta_{\pi}(g_1,g_2)^{-1}\cdot
\la \pi(g_1^{-1}),\lambda{(g_2)}\ra \cdot\zeta'_{\pi}(g_1,g_2)
$$
is inflated from $Q$.\\
(4) Suppose $\pi\in B^1(Q,\ach)$. Then there is a $\chi_{\pi}\in
\ach$ such that $\pi(g)=g(\chi_{\pi})\chi_{\pi}^{-1}$ for every
$g\in Q$. Setting
$\zeta_{\pi}(g_1,g_2):=\la \chi_{\pi},\beta(g_1,g_2)\ra $, we obtain
$$
\la \pi(g_1),g_1(\beta(g_2,g_3))\ra =\la g_1(\chi_{\pi})\chi_{\pi}^{-1},g_1(\beta(g_2,g_3))\ra
=\delta_Q^2(\zeta_{\pi})(g_1,g_2,g_3).
$$
It follows that
\begin{equation*}
\begin{split}
c_{\pi}(a_1\bar{g}_1,a_2\bar{g}_2)
=&~\varphi_{\pi}
\bar{\zeta}_{\pi}(a_1\bar{g}_1,a_2\bar{g}_2)=
\langle\pi(g_1),g_1(a_2)\rangle^{-1}\la \chi_{\pi},\beta(g_1,g_2)\ra \\
=&~\la g_1(\chi_{\pi}),g_1(a_2)\ra ^{-1}\la \chi_{\pi},g_1(a_2)\ra \la \chi_{\pi},\beta(g_1,g_2)\ra \\
=&~\la \chi_{\pi},a_1a_2\ra ^{-1}\la \chi_{\pi},a_1g_1(a_2)\beta(g_1,g_2)\ra
=\delta_G^1(\bar{\chi}^{-1}_{\pi})(a_1\bar{g}_1,a_2\bar{g}_2),
\end{split}
\end{equation*}
where $\bar{\chi}_{\pi}\in C^1(G,\C^*)$ is defined by
$\bar{\chi}_{\pi}(a\bar{g})=\chi_{\pi}(a)$. Hence, with the above
choice of $\zeta_{\pi}$, we deduce that $c_{\pi}\in B^2(G,\C^*)$.\\
(5) This follows from the definition of $c_\pi$ (\ref{cp}) and from equation (\ref{varphimulti}).
\end{proof}

We summarize the above in:
\begin{corollary}\label{C}
The map $C:[\pi]\mapsto [c_{\pi}]\mod[\im(\inf^Q_G)]$ is a well
defined homomorphism from $\kb$ to
$\ker(\res^G_A)/\im(\inf^Q_G)$.
\end{corollary}

\subsection{}\label{subsecPi}
We show that $C$ is actually an isomorphism by presenting its inverse map
$\Pi:\ker(\res^G_A)/\im(\inf^Q_G)\to \kb$.
A similar construction
can be found in \cite[\S 1.7]{K}.

\begin{lemma}\label{satkarp}
With the above notation, any class in $\ker(\res^G_A)$ admits a representative $c\in
Z^2(G,\C^*)$ such that
\begin{equation}\label{karp}
c(\gamma_{1},\gamma_{2})=
c(a \gamma_{1},\gamma_{2})
\end{equation}
for all $\gamma_{1},\gamma_{2}\in G$, and $a \in A$.
\end{lemma}

\begin{proof}
See \cite[Lemma 7.1, P.59]{K}.
\end{proof}

Let $c\in Z^2(G,\C^*)$ be a 2-cocycle on $G$ such that
$[c]\in\ker(\res^G_A)$. We need to construct a 1-cocycle on $Q$
with values in $\ach$. For this purpose we first show how $[c]\in
\ker(\res^G_A)$ yields a character on $A$ for any $g\in Q$.
By Lemma \ref{satkarp}, we may assume
 $c(\gamma_{1},\gamma_{2})={c(a\gamma_{1},\gamma_{2})}$ for all $\gamma_{1},\gamma_{2}
 \in G$ and $a\in A$.
Now, for any $a\in A$ and any $g\in Q$ define
\begin{eqnarray}\label{pi}
\begin{array}{rcl}
\pi_c(g):A & \to & \C^*\\
a & \mapsto & c({\overline{g^{-1}},a}).
\end{array}
\end{eqnarray}

\begin{proposition}\label{prop}
Let $\pi_c(g)$ be as in (\ref{pi}). Then:
\begin{enumerate}
\item For any $g\in Q$, $\pi_c(g)\in \ach$.
\item $\pi_c(g)$ does not depend on the choice of the transversal set $\{\bar{g}\}_{g\in Q}$.
\item  The map $g\mapsto \pi_c(g)$ is a 1-cocycle from $Q$ to $\ach$.
\item If $c'\in [c]$ satisfies (\ref{karp}),
then $\pi_{c'}\pi_c^{-1}$
is a coboundary.
\item If $c_1$ and $c_2$ satisfy (\ref{karp}), then so does $c_1c_2$. Moreover, $\pi_{c_1c_2}=\pi_{c_1}\pi_{c_2}$.
\item $\pi_c\in B^1(Q,\ach)$ if and only if $[c]\in \im(\inf^Q_G)$.
\item
The 3-cochain $(g_1,g_2,g_3)\mapsto \la \pi_c(g_1),g_1(\beta(g_2,g_3))\ra $
is a coboundary on $Q$.
\end{enumerate}
\end{proposition}
\begin{proof}
See \cite[Theorem 7.3, P.60]{K} for the proof of (1)-(6).
Let us prove (7).
By Proposition \ref{vp}, $\la \pi(g_1),g_1(\beta(g_2,g_3))\ra =\delta_G^2(\varphi_{\pi}^{-1})
(a_1\bar{g}_1,a_2\bar{g}_2,a_3\bar{g}_3)$
for any $\pi\in Z^1(Q,\ach)$ and $a_1\bar{g}_1,a_2\bar{g}_2,a_3\bar{g}_3\in G$. We
claim that in case $\pi=\pi_c$, we can substitute the $G$-cochain
$\varphi_{\pi_c}$ with an appropriate $Q$-cochain.
Indeed,
\begin{equation}\label{neweq1}
\varphi_{\pi_c}(a_1\bar{g}_1,a_2\bar{g}_2)=\la \pi_c(g_1),g_1(a_2)\ra ^{-1}=
\la \pi_c(g_1^{-1}),a_2)\ra =c(\bar{g_1},a_2).
\end{equation}
Applying the 2-cocycle condition, we obtain
\begin{eqnarray}\label{neweq2}
\begin{array}{l}
c(\bar{g_1},a_2)=
c(\bar{g_1},a_2\bar{g_2})c(\bar{g_1}a_2,\bar{g_2})^{-1}c(a_2,\bar{g_2})=\\
=c(\bar{g_1},a_2\bar{g_2})c(g_1(a_2)\bar{g_1},\bar{g_2})^{-1}c(a_2,\bar{g_2}).
\end{array}
\end{eqnarray}
Since $c$ satisfies (\ref{karp}), then equations (\ref{neweq1}) and (\ref{neweq2}) imply
\begin{equation}\label{deltaeq}
\varphi_{\pi_c}(a_1\bar{g}_1,a_2\bar{g}_2)=
c(a_1\bar{g_1},a_2\bar{g_2})c(\bar{g_1},\bar{g_2})^{-1}.
\end{equation}
Note that since $c$ is a
2-cocycle on $G$, we have $\delta_G^2(c)=0$. From (\ref{deltaeq}) we obtain that
$$\delta_G^2(\varphi_{\pi_c})(a_1\bar{g}_1,a_2\bar{g}_2,a_3\bar{g}_3)
=\delta_Q^2(\tilde{c}^{-1})(g_1,g_2,g_3),$$ where $\tilde{c}\in C^2(Q,\C^*)$ is defined by
$\tilde{c}(g,g'):=c(\bar{g},\bar{g}'),\ \ g,g'\in Q$.
In particular,
$$\la \pi_c(g_1),g_1(\beta(g_2,g_3))\ra =\delta_Q^2(\tilde{c})(g_1,g_2,g_3).$$
\end{proof}
We summarize the above.
\begin{corollary} \label{Pi}
The map $\Pi:[c]\!\!\mod[\im(\inf^Q_G)] \mapsto [\pi_c]$ is a well
defined homomorphism from $\ker(\res^G_A)/\im(\inf^Q_G)$ to
$\kb\subseteq H^1(Q,\ach)$.
\end{corollary}
To complete the discussion we have:
\begin{proposition}
The homomorphisms $C$ and $\Pi$ in Corollaries \ref{C} and \ref{Pi} are mutually inverse.
\end{proposition}
\begin{proof}
Let $[\pi]\in \kb$.
Note that since $\beta$ is normalized, then by (\ref{vs}) we have
$$\delta_Q^2(\zeta_{\pi})(g,1,g)=1$$
for every $g\in Q$.
Developing $\delta_Q^2$ we obtain that
\begin{equation}\label{1gg1}
\zeta_{\pi}(1,g)=\zeta_{\pi}(g,1), \forall g\in Q.
\end{equation}
We first claim that the choice of the cochain $\zeta_{\pi}$ satisfying (\ref{vs}) can be done such that
\begin{equation}\label{normvs}
 \zeta_{\pi}(1,g)=\zeta_{\pi}(g,1)=1, \forall g\in Q.
\end{equation}
Indeed, let
$$\zeta'_{\pi}(g_1,g_2):=\zeta_{\pi}(g_1,g_2)\zeta_{\pi}(1,g_2)^{-1}.$$
Then clearly $$\zeta'_{\pi}(1,g)=1$$ for every $g\in Q$.
To complete the proof of the claim, we need to check that $\zeta'_{\pi}$ satisfies equation (\ref{vs}).
Indeed, for every $g_1,g_2,g_3\in Q$,
$$\begin{array}{l}
\delta_Q^2(\zeta'_{\pi})(g_1,g_2,g_3)=\delta_Q^2(\zeta_{\pi})(g_1,g_2,g_3)
\zeta_{\pi}^{-1}(1,g_2g_3)\zeta_{\pi}(1,g_2)=\\
\delta_Q^2(\zeta_{\pi})(g_1,g_2,g_3)
\delta_Q^2(\zeta_{\pi}^{-1})(1,g_2,g_3)
=\la \pi(g_1),g_1(\beta(g_2,g_3))\ra \la \pi(1),\beta(g_2,g_3)\ra ^{-1}=\\
=\la \pi(g_1),g_1(\beta(g_2,g_3))\ra (=\delta_Q^2(\zeta_{\pi})(g_1,g_2,g_3)).
\end{array}$$
Note that the 2-cocycle $c_{\pi}:=\varphi_{\pi}\bar{\zeta}_{\pi}$ satisfies condition (\ref{karp}),
and hence $\pi_{c_{\pi}}$ is well defined.
By the definitions (\ref{obt}),(\ref{cp}) and (\ref{pi}), we have that
for every $g\in Q$ and $a\in A$,
$$\la \pi_{c_{\pi}}(g),a\ra =c_{\pi}(\overline{g^{-1}},a)=\la \pi(g),a\ra \cdot\zeta_{\pi}(g^{-1},1).$$
Choosing $\zeta_{\pi}$ that satisfies (\ref{normvs}), we obtain $\la \pi_{c_{\pi}}(g),a\ra =\la \pi(g),a\ra $.
By Proposition \ref{condpi}(2) and Proposition \ref{prop}(6) we deduce that
$$[\pi_{c_{\pi}}]=[\pi].$$

Conversely, let $c\in Z^2(G,\C^*)$ be as in Lemma \ref{satkarp}.
Then for every $g_1,g_2\in Q$ and $a_1,a_2\in A$,
$$c_{\pi_c}(a_1\bar{g}_1,a_2\bar{g}_2)=\varphi_{\pi_c}(a_1\bar{g}_1,a_2\bar{g}_2)\zeta_{\pi_c}(g_{1},g_2).$$
By equation (\ref{deltaeq}),
$$c_{\pi_c}(a_1\bar{g}_1,a_2\bar{g}_2)=
c(a_1\bar{g_1},a_2\bar{g_2})c(\bar{g_1},\bar{g_2})^{-1}\zeta_{\pi_c}(g_{1},g_2).$$
This shows that $c$ and $c_{\pi_c}$ differ by a 2-cocycle which is inflated from
$Q$. Thus, $$[c_{\pi_c}]=[c]\mod[\im(\infl^Q_G)].$$
\end{proof}

\subsection{}
We can interpret the 1-cocycle $\pi_c$ as interchanging
the primitive idempotents of $\C[A]$ by $Q$ via the conjugation in $\C^c[G]$ as in
Proposition \ref{pinterchange} below. This will be needed in the next section. First, some notation:

Let $c\in Z^2(G,\C^*)$ be any 2-cocycle and let $\C^c[G]$
be the corresponding twisted group algebra with a $\C$-basis
$\{U_{\gamma}\}_{\gamma\in G}$
(satisfying $U_{\sigma}U_{\gamma}=c(\sigma,\gamma)U_{\sigma\gamma}$).
Then $\gamma\in G$ acts on $\C^{c}[G]$ by left conjugation with $U_{\gamma}$:
\begin{equation}\label{action}
\gamma \cdot U_{\sigma}=U_{\gamma}U_{\sigma}U_{\gamma}^{-1}=[\sigma,\gamma]_cU_{\gamma\sigma\gamma^{-1}},
\end{equation}
where
\begin{equation}\label{alph}
[\sigma,\gamma]_c:=c(\gamma,\sigma)c(\gamma\sigma\gamma^{-1},\gamma)^{-1}
=c(\gamma,\sigma\gamma^{-1})c(\sigma\gamma^{-1},\gamma)^{-1} .
\end{equation}
Next, let $[\beta]: 1\to A\to G\to Q\to 1$ be an extension of finite groups,
and let $c\in Z^2(G,\C^*)$ be any 2-cocycle.
Then $G$ acts on the central primitive idempotents of the subalgebra
$\C^{\res^G_A(c)}[A]\subseteq \C^c[G]$ by conjugation.
For every $\gamma\in G$ and a central primitive idempotent $\iota\in\C^{\res^G_A(c)}[A]$,
$$\gamma(\iota)=U_{\gamma}\cdot\iota \cdot U_{\gamma}^{-1}.$$

Suppose now that $A$ is abelian and $[c]\in \ker(\res^G_A)$.
By Lemma \ref{satkarp} we may assume that $c$ satisfies (\ref{karp}).
In this case the commutative subalgebra $\C[A]=\C^{\res^G_A(c)}[A]\subseteq \C^c[G]$
admits the primitive idempotents
$$\iota_{\chi}=\frac{1}{|A|}\sum\limits_{a\in A}\la \chi,a\ra ^{-1}U_a$$
for every ${\chi\in \ach}$.
The action of $G$ on the set of primitive idempotents is via the quotient $Q$.
Then the 1-cocycle $\pi_c$ associated to $c$ determines this action as follows.

\begin{proposition}\label{pinterchange}
With the above notation let $g\in Q$ and $\chi\in \ach$.
Then
\begin{equation}\label{charact}
{g}(\iota_{\chi})=\iota_{g(\chi)\pi_c(g)}.
\end{equation}
\end{proposition}

\begin{proof}
\begin{equation*}
\begin{split}
{g}(\iota_{\chi})=&~U_{\bar{g}}\iota_{\chi}U_{\bar{g}}^{-1}=\frac {1}{|A|}
\sum\limits_{a\in  A}\la \chi,a\ra ^{-1}U_{\bar{g}}U_a U_{\bar{g}}^{-1}\\
=&~\frac {1}{|A|} \sum\limits_{a\in  A}\la \chi,a\ra ^{-1}[a,\bar{g}]_cU_{g(a)}
=~\frac {1}{|A|} \sum\limits_{a\in  A}\la \chi,a\ra ^{-1}\la \pi_c(g^{-1}),a\ra \cdot U_{g(a)}\\
=&~ \frac {1}{|A|} \sum\limits_{a\in  A}\la g(\chi)\pi_c(g),g(a)\ra ^{-1} U_{g(a)}
=\iota_{g(\chi)\pi_c(g)}.
\end{split}
\end{equation*}
\end{proof}

\begin{remark}\label{genpi}
This construction of the map $\Pi$ in \S \ref{subsecPi} works also if
$\C^*$ is replaced by any $A$-trivial $G$-module.
In particular, assume $G$ as above acts on a field $K$, with $A$ acting trivially, and
let $c\in Z^2(G,K^*)$ be a 2-cocycle such that $[c]\in \ker(\res^G_A)$.
One can then establish (\ref{charact}) for the action of $Q$ on the primitive
idempotents of $K[A]$ in the crossed product $K^c*[G]$.
\end{remark}

\section{Bijective and Non-Degenerate Cohomology Classes} \label{2main}

In this section we prove:
\begin{enumerate}
\item If $\pi:\:Q \rightarrow \ach$ is a bijective 1-cocycle such that
$[\pi]\in \kb$, then $c_{\pi}\in Z^{2}(G,{\C}^{*})$ is a non-degenerate 2-cocycle
(Theorem \ref{P:122}).
\item If $c\in Z^2(G,\C^*)$ is a non-degenerate 2-cocycle satisfying (\ref{karp})
(in particular $[c]\in \ker(\res^G_A)$) and if $|A|=|Q|$, then $\pi_c$ is bijective
(Theorem \ref{thm:pic_bi}).
\end{enumerate}

This will establish the correspondence in Theorem A between bijective classes and non-degenerate
classes $\mod[\im(\inf^Q_G)]$, once we have shown that at least one of the two properties
is independent of the choice of cocycle. We therefore prove:

\begin{proposition}
Let $\pi\in Z^1(Q,M)$, where $M$ is any $Q$-module.
If $\pi$ is surjective, then so are all the cocycles in its cohomology class.
\end{proposition}

\proof
Let $\pi'\in Z^1(Q,M)$ such that $[\pi']=[\pi]$.
Then there exists $m\in M$ such that for every $g\in Q$, $\pi'(g)=\pi(g)g(m)m^{-1}$.
By the hypothesis, $m=\pi(g')$ for some $g'\in Q$.
Then $$\pi'(g)=\pi(g)g(\pi(g'))m^{-1}=\pi(gg')m^{-1}.$$ Hence $\pi'$ is surjective as well.
\qed\\

Let $[\beta]: 1\to A\to G\to Q\to 1$ be an extension of finite groups,
where $A$ is abelian.
We first prove that if $\pi:Q \rightarrow \ach$ is a bijective 1-cocycle
(in particular $|Q|=|A|$) such that
$[\pi]\in \kb$, then the twisted group algebra
$\C^{c_{\pi}}[G]$ is isomorphic as a $G$-module
by conjugation (\ref{action}) to the regular left representation $\C[G]$ of $G$.
As mentioned in the introduction, this will show that $c_{\pi}\in Z^{2}(G ,{\C}^{*})$ is a
non-degenerate 2-cocycle [DMJ, Corollary 3].

Let $c\in Z^2(G,\C^*)$ be any 2-cocycle and let
$\{U_{\gamma}\}_{\gamma\in G}$ be a $\C$-basis of the left $G$-module $\C^c[G]$,
where $\gamma\in G$ acts by conjugation (\ref{action}).
Further, let $\{\gamma\}_{\gamma\in G}$ be a $\C$-basis of the left $G$-module $\C[G]$,
where $\gamma\in G$ acts by left multiplication.
We present two morphisms between these $G$-modules, which are mutually inverse
in case $c=c_{\pi}$
for a bijective 1-cocycle $\pi$ such that $[\pi]\in \kb$
and hence yield the $G$-module isomorphism between $\C^{c_{\pi}}[G]$ and $\C[G]$.

With the notation of (\ref{alph}) we define:
 \begin{eqnarray} \label{psi}
 \begin{array} {rcl}
 \psi_c:\C^{c}[G]& \rightarrow & \C[G]         \\
  U_{\gamma}& \mapsto &  \sum\limits_{\sigma\in G }[\gamma,\sigma^{-1}]_c\sigma
\end{array}
\end{eqnarray}
and
\begin{eqnarray} \label{theta}
\begin{array} {rcl}
\theta_c: \C[G]&\rightarrow & \C^{c}[G] \\
\sigma& \mapsto & \frac{1}{|G |}\sum\limits_{\gamma\in G }[\gamma,\sigma^{-1}]_c^{-1}U_{\gamma}
\end{array}
\end{eqnarray}
                (and extend
                $\psi_c$ and $\theta_c$ linearly).

        \begin{proposition}
            The maps $\psi_c$ and $\theta_c$ are $G$-module morphisms.
        \end{proposition}
        \proof For any $\gamma,\tau\in G$ we have
$$\begin{array}{l}
\psi_c(\tau\cdot U_{\gamma})=[\gamma,\tau]_c\psi_c(U_{\tau\gamma \tau^{-1}})=
[\gamma,\tau]_c\sum\limits_{\sigma\in \Gamma }[\tau\gamma \tau^{-1},\sigma^{-1}]_c\sigma=\\
\sum\limits_{\sigma\in
                \Gamma }[\gamma,\sigma^{-1}\tau]_c\sigma = \sum\limits_{\eta\in
                \Gamma }[\gamma,\eta^{-1}]_c\tau \eta=\tau\sum\limits_{\eta\in \Gamma }
                [\gamma,\eta^{-1}]_c\eta=\tau\psi_c (U_{\gamma}).
               \end{array}$$
           Hence $\psi_c$ is a $G$-module map.

            Next,
$$ \theta_c(\tau)=\frac{1}{|G |}\sum\limits_{\gamma\in G }[\gamma,\tau^{-1}]_c^{-1}U_{\gamma}=
\frac{1}{|G|}\sum\limits_{\gamma\in G }[\tau\gamma\tau^{-1},\tau^{-1}]_c^{-1}
U_{\tau\gamma\tau^{-1}}=$$
by definition (\ref{alph})
$$\begin{array}{l}
=\frac{1}{|G|}\sum\limits_{\gamma\in G }c(\tau^{-1},\tau\gamma\tau^{-1})^{-1}
c(\gamma,\tau^{-1})U_{\tau\gamma\tau^{-1}}=\\
\frac{1}{|G|}\sum\limits_{\gamma\in G }c(1,\gamma\tau^{-1})^{-1}
c(\tau,\gamma\tau^{-1})c(\tau^{-1},\tau)^{-1}c(\gamma,\tau^{-1})U_{\tau\gamma\tau^{-1}}=\\
\frac{1}{|G|}\sum\limits_{\gamma\in G }
c(\tau,\gamma\tau^{-1})c(\tau^{-1},\tau)^{-1}c(\gamma,\tau^{-1})c(\gamma,1)^{-1}
U_{\tau\gamma\tau^{-1}}=\\
\frac{1}{|G|}\sum\limits_{\gamma\in G }
c(\tau,\gamma\tau^{-1})c(\gamma\tau^{-1},\tau)^{-1}U_{\tau\gamma\tau^{-1}}=
\end{array}$$
again, by (\ref{alph})
$$=\frac{1}{|G |}\sum\limits_{\gamma\in G }[\gamma,\tau]_cU_{\tau\gamma\tau^{-1}}
=\frac{1}{|G|}\sum\limits_{\gamma\in
                G }\tau\cdot U_{\gamma}= \tau\cdot\theta_c(1).$$
            Thus,
            $$\theta_c(\sigma\tau)=\sigma\tau\cdot \theta_c(1)=\sigma\cdot
(\tau\cdot\theta_c(1))=\sigma\cdot\theta_c(\tau)$$ and $\theta_c$ is a $G$-module map.
This completes the proof of the proposition.\qed\\

We now give a sufficient condition on the 2-cocycle $c$ for
$\psi_c$ and $\theta_c$ to be inverse to each other.

    \begin{proposition}\label{T:121}
Let $c\in Z^{2}(G ,\C^{*})$ be such that for any non-trivial $\tau\in G $,
$\sum\limits_{\gamma\in G }[\gamma,\tau]_c=0$.
Then $\psi_c$ is invertible with $\psi_c^{-1}=\theta_c$.
In particular, $c$ is non-degenerate.
    \end{proposition}
    \proof
            Since dim$_{\C}\C^{c}[G]$ = dim$_{\C}\C[G]$ =
            $|G|$, it is enough to show that $\theta_c$
            is a right inverse of $\psi_c$. We first compute $\psi_c\circ \theta_c$ on
            the trivial element $1\in G$.
\begin{equation}
                (\psi_c\circ \theta_c)(1)=\psi_c(\frac{1}{|G |}\sum_{\gamma\in G }
                U_{\gamma})=\frac{1}{|G|}\sum_{\gamma\in G }\sum_{\sigma\in
                G }[\gamma,\sigma^{-1}]_c\sigma\,.
\end{equation}
By the hypothesis,
$\frac{1}{|G |}\sum\limits_{\gamma\in G}[\gamma,\sigma^{-1}]_c=\delta_{1,\sigma}$.
Thus $$(\psi_c\circ\theta_c)(1)=\sum\limits_{\sigma\in G }\delta_{1,\sigma}\sigma=1.$$
Consequently,
$$(\psi_c\circ\theta_c)(\tau)=\psi_c(\theta_c(\tau))=\psi_c(\tau\cdot\theta_c(1))=
\tau\psi_c(\theta_c(1))=\tau \cdot 1=\tau.$$\qed

  \begin{theorem}\label{P:122}
Let $\pi:Q \to \ach$ be a bijective 1-cocycle such that
$[\pi]\in \kb$,
            then $c_{\pi}\in Z^{2}(G,{\C}^{*})$ is a non-degenerate 2-cocycle.
    \end{theorem}
    \proof
We show that $c_{\pi}\in Z^2(G,\ach)$ satisfies the criterion in Proposition \ref{T:121}.
Recall that for every $a_1\bar{g}_1,a_2\bar{g}_2\in G$,
$$c_{\pi}(a_1\bar{g}_1,a_2\bar{g}_2)=\la \pi(g_1^{-1}),a_2\ra \zeta_{\pi}(g_1,g_2),$$
where $\zeta_{\pi}\in C^2(Q,\C^*)$ satisfies (\ref{vs}).
Let $a\bar{g}$ be a non-trivial element in G. We prove that
$\sum_{\gamma\in G}[\gamma,a\bar{g}]_{c_{\pi}}=0$.
 $$\begin{array}{l}
  \sum\limits_{a'\bar{g}'\in G }[a'\bar{g}',a\bar{g}]_{c_{\pi}}=
  \sum\limits_{a'\bar{g}'\in G }c_{\pi}(a\bar{g},a'\bar{g}')
  c_{\pi}((a\bar{g})a'\bar{g}'(a\bar{g})^{-1},a\bar{g})^{-1}=\\
\sum\limits_{a'g'\in G }\la \pi(g^{-1}),a'\ra \zeta_{\pi}(g,g')
\la \pi(gg'^{-1}g^{-1}),a\ra ^{-1}\zeta_{\pi}(gg'g^{-1},g)^{-1}=\\
(\sum\limits_{a'\in A}\la \pi(g^{-1}),a'\ra )\cdot(\sum\limits_{g'\in Q}\zeta_{\pi}(g,g')
\la \pi(gg'^{-1}g^{-1}),a\ra ^{-1}\zeta_{\pi}(gg'g^{-1},g)^{-1}).
\end{array}$$
If $g$ is not trivial in $Q$ then the injectivity of $\pi$ implies that
$\pi(g^{-1})$ is not the trivial character in $\ach$.
Consequently
\begin{equation}\label{gnontriv}
\sum_{a'\in A}\la \pi(g^{-1}),a'\ra =0, \ \ g\neq 1.
\end{equation}
Suppose then that $g=1$. Applying (\ref{1gg1}) we obtain
\begin{eqnarray}
\begin{array}{l}
\sum_{g'\in Q}\zeta_{\pi}(g,g')
\la \pi(gg'^{-1}g^{-1}),a\ra ^{-1}\zeta_{\pi}(gg'g^{-1},g)^{-1}=\\
\sum_{g'\in Q}\la \pi(g'^{-1}),a\ra ^{-1}.
\end{array}
\end{eqnarray}
Since $\pi:Q \to \ach$ is bijective, we have
\begin{equation}\label{gtriv}
\sum_{g'\in Q}\la \pi(g'^{-1}),a\ra ^{-1}=\sum_{\chi\in \ach}\la \chi,a\ra ^{-1}=|A|\delta_{a,1}.
\end{equation}
From (\ref{gnontriv}) and (\ref{gtriv}) we see that
$\sum\limits_{\gamma\in G }[\gamma,a\bar{g}]_{c_{\pi}}$ vanishes unless both
$g$ and $a$ are trivial.
\qed\\

We now prove the other direction of the correspondence in Theorem A.

Let $[\beta]: 1\to A\to G\to Q\to 1$ be an extension of finite groups,
and let $c\in Z^2(G,\C^*)$ be any 2-cocycle.
Let $\iota$ be any central primitive idempotent of the subalgebra
$\C^{\res^G_A(c)}[A]\subseteq\C^c[G]$. It is easily checked that
$\C^c[G](\bigoplus\limits_{\gamma\in G}\gamma(\iota))$ is a two sided ideal of $\C^c[G]$.
We obtain a necessary condition for the simplicity of $\C^c[G]$, or equivalently,
for the non-degeneracy of $c$.
See [AS, Lemma 1.2] for a more general setup.
\begin{proposition}\label{transitive}
If a 2-cocycle $c\in Z^2(G,\C^*)$ is non-degenerate, then $G$ acts
transitively on the set of central primitive idempotents of the subalgebra
$\C^{\res^G_A(c)}[A]$ for any normal subgroup $A\lhd G$.
\end{proposition}

Assume again that $A$ is abelian and $[c]\in \ker(\res^G_A)$ satisfies (\ref{karp}).
Substituting
$\iota_{1}=\frac {1}{|A|} \sum\limits_{a\in  A}U_a$, the primitive idempotent of $\C[A]$
which corresponds to the trivial character $1_A\in \ach$
into equation (\ref{charact}) we obtain
\begin{equation}\label{iot}
{g}(\iota_{1})=\iota_{\pi_c(g)}.
\end{equation}
Finally, equation (\ref{iot}) and Proposition \ref{transitive} imply:
\begin{theorem}\label{thm:pic_bi}
Let $1\to A\to G\to Q\to 1$ be an extension of finite groups with $A$ abelian
and let $c\in Z^2(G,\C^*)$ be a non-degenerate 2-cocycle satisfying (\ref{karp}).
Then $\pi_c\in Z^1(Q,\ach)$ is surjective.
In particular, if $|A|=|Q|$, then $\pi_c$ is bijective.
\end{theorem}
This completes the proof of Theorem A.

\section{Example}\label{Example}
The following example was shown to us by E. Aljadeff and \'{A}. del R\'{\i}o.\\
Let $Q:=C_4\times C_2=\la \sigma\ra \times\la \tau\ra $ act on $A:=C_4\times C_2=\la x\ra \times\la y\ra $ as follows:
\begin{equation}
\sigma(x)=xy,\ \ \tau(x)=x^{-1},\ \ \sigma(y)=\tau(y)=y.
\end{equation}

Let $\ach:=\la \check{x}\ra \times\la \check{y}\ra $, where
$\check{x}:\begin{array}{rcr}x & \mapsto & i \\ y & \mapsto & 1\end{array}$ and
$\check{y}:\begin{array}{rcr}x & \mapsto & 1\\ y & \mapsto & -1\end{array}$.\\
Define:
\begin{eqnarray}
\begin{array}{rcl}
\pi:Q & \to & \ach \\
\sigma^k\tau^l & \mapsto & \check{x}^{l-(-1)^lk}\check{y}^l.
\end{array}
\end{eqnarray}
One can easily check that $\pi$ is a bijective $1$-cocycle.

Next, consider the extension:
$$[\beta]:\; 1\to A\to G:=\la x,y,\bar{\sigma},\bar{\tau}\ra \to Q\to 1$$
determined by
$$\bar{\sigma}^4=\bar{\tau}^2=1,\ \
\bar{\tau}\bar{\sigma}\bar{\tau}^{-1}\bar{\sigma}^{-1}=x^{-1}y.$$
Then $[\beta] \cup [\pi]=0$, and by Corollary B, the group $G$
is of central type. However, there is no abelian normal subgroup
$N\lhd G$ of order $8$ such that $G:=N\rtimes (G/N)$. Corollary B
therefore allows us to construct central type groups which
could not be obtained by [EG3].

\section{Twists for Finite Groups}\label{twist}
Theorem A has a useful application to the explicit
construction of twists for finite groups.
Let $\phi:\C[G]\to C'$ be a $G$-module isomorphism between the regular $G$-module $\C[G]$
and a $G$-coalgebra $(C',\bigtriangleup',\varepsilon')$.
A result of Movshev [M, Proposition 5] says that
$C'$ is isomorphic as a $G$-coalgebra to the coalgebra $\C[G]_{J_{\phi}}$, where
\begin{equation}\label{mov}
J_{\phi}=(\phi^{-1}\otimes\phi^{-1})\bigtriangleup'(\phi(1))
\end{equation}
is a twist for $G$.

Let $[\beta]: 1\to A\to G\to Q\to 1$ be an extension, and suppose that
$[\pi]\in H^1(Q,\ach)$ is a bijective class such that $[\beta] \cup [\pi]=0$.

Let $c_{\pi}$ be a 2-cocycle on $G$ determined by $\pi$ as in equation (\ref{cp}) above.
We apply (\ref{mov}) to the $G$-coalgebra $C'=(\C^{c_{\pi}}[G])^*$ by a $G$-module isomorphism
$\phi_{\pi}:\C[G]\to (\C^{c_{\pi}}[G])^*$ described as follows.
Equation (\ref{theta}) above yields
an explicit $G$-module homomorphism $\theta_{c_{\pi}}$
between the regular module $\C[G]$ and $\C^{c_{\pi}}[G]$ with the conjugation action.
Under the above assumption on $\pi$, the proof of Theorem \ref{P:122}
(claiming that $c_{\pi}$ is non-degenerate) shows that
$\theta_{c_{\pi}}$ is actually an isomorphism.
The element $r_{\pi}=\theta_{c_{\pi}}(1)\in \C^{c_{\pi}}[G]$ is {\it regular},
that is $\{\gamma(r_{\pi})\}_{\gamma\in G}$ is a $\C$-basis of $\C^{c_{\pi}}[G]$.
Similarly, the element $r^*_{\pi}\in (\C^{c_{\pi}}[G])^*$ defined by:
\begin{equation}
r^*_{\pi}(\gamma(r_{\pi}))=\delta_{\gamma,e}
\end{equation}
generates a $\C$-basis of $(\C^{c_{\pi}}[G])^*$.
Then
\begin{eqnarray}
\begin{array}{rcl}
\phi_{\pi}:\C[G] & \to & (\C^{c_{\pi}}[G])^*\\
\gamma& \mapsto & \gamma(r^*_{\pi})
\end{array}
\end{eqnarray}
is a $G$-module isomorphism and $J_{\phi_{\pi}}$ is a twist for $G$.

\end{document}